# A study on radial basis function and quasi-Monte Carlo methods

W. Chen[*] and J. He[**]

[*]Department of Mechanical System Engineering, Shinshu University, Wakasato 4-17-1, Nagano City, Nagano 380-8533, Japan (E-mail: chenw@homer.shinshu-u.ac.jp)

[**]P.O. Box 189, Shanghai University, 149 Yanchang Road, Shanghai 200072, P.R. China

## 1. Introduction

The radial basis function (RBF) and quasi Monte Carlo (QMC) methods are two very promising schemes to handle high-dimension problems with complex and moving boundary geometry due to the fact that they are independent of dimensionality and inherently meshless. The two strategies are seemingly irrelevant and are so far developed independently. The former is largely used to solve partial differential equations (PDE), neural network, geometry generation, scattered data processing with mathematical justifications of interpolation theory [1], while the latter is often employed to evaluate high-dimension integration with the Monte Carlo method (MCM) background [2].

The purpose of this communication is to try to establish their intrinsic relationship on the grounds of numerical integral. The kernel function of integral equation is found the key to construct efficient RBFs. Some significant results on RBF construction, error bound and node placement are also presented. It is stressed that the RBF is here established on integral analysis rather than on the sophisticated interpolation and native space analysis.

## 2. Relationship between RBF and QMC

In this section, we try to relate the RBF and QMC through evaluation of particular solution of PDE. Without loss of generality, consider Possion equation

$$\nabla^2 u = f(x), x \in \Omega, \quad (1)$$

$$u(x) = D(x), x \subset S_u, \quad (2a)$$

$$\frac{\partial u(x)}{\partial n} = N(x), x \subset S_T, \quad (2b)$$

where $x$ means multi-dimensional independent variable. $n$ is the unit outward normal. The solution of Eq. (1) can be expressed as

$$u = v + u_p, \quad (3)$$

where $v$ and $u_p$ are the general and particular solutions, respectively. The latter satisfies

$$\nabla^2 u_p = f(x) \quad (4)$$

but does not necessarily satisfy boundary conditions. $v$ is the homogeneous solution

$$\nabla^2 v = 0, \quad (5)$$

$$v(x) = D(x) - u_p, \quad (6a)$$

$$\frac{\partial v(x)}{\partial n} = N(x) - \frac{\partial u_p(x)}{\partial n}, \quad (6b)$$

Eqs. (5) and (6a,b) can be solved by various boundary numerical techniques such as the BEM, method of fundamental solution (MFS) [3] and boundary knot method (BKM) [4]. In this study, however, we pay attention only on the evaluation of the particular solution.

As was pointed out by Chen et al. [3], the dual reciprocity method (DRM) and QMC are two general techniques evaluating particular solution. What follows is a brief description of these two strategies with emphasis on the analysis of the RBF and MQC.

### 2.1. DRM, RBF and integral equation

In terms of the DRM, the inhomogeneous terms of Eq. (4) are approximated at first by

$$f(x) \cong \sum_{k=1}^{N+L} \alpha_k \phi_k(r_k) \qquad (7)$$

where $\alpha_k$ are the unknown coefficients. $N$ and $L$ are respectively the numbers of knots on domain and boundary. $r_k = \|x - x_k\|$ represents the Euclidean distance norm, $\phi$ is the RBF and satisfies

$$\phi(r) = \nabla^2 \psi(r). \qquad (8)$$

We can uniquely determine

$$\alpha = A_\phi^{-1} f(x), \qquad (9)$$

where $A_\phi$ is nonsingular RBF interpolation matrix. Finally, we can get particular solutions at any point by summing localized approximate particular solutions

$$u_p = \sum_{k=1}^{N+L} \alpha_k \psi_k(r_k). \qquad (10)$$

where $\phi$ is also RBF. $u_p$ is evaluated by substituting Eq. (9) into Eq. (10).

Another DRM-equivalent method is to compute the domain integral. By using Green second theorem, the solution of Eq. (1) can be expressed as

$$u(x) = \int_\Omega f(z) u^*(x,z) d\Omega(z) + \int_\Gamma \left\{ u \frac{\partial u^*(x,z)}{\partial n(z)} - \frac{\partial u}{\partial n(z)} u^*(x,z) \right\} d\Gamma(z), \qquad (11)$$

where $u^*$ is the fundamental solution of Laplace operator. $z$ denotes source point. It is noted that the first and second terms of Eq. (11) respectively equal the particular and general solutions of Eq. (3). If a numerical integral scheme is used to approximate Eq. (11), we have

$$u(x) \cong \sum_{k=1}^{N+L} \omega(x, x_k) f(x_k) u^* + \sum_{k=N+1}^{N+L} \omega(x, x_k) \left[ u \frac{\partial u^*}{\partial n} - \frac{\partial u}{\partial n} u^* \right], \qquad (12)$$

where $\omega(x, x_j)$ are the weights dependent on the integral schemes. Furthermore, we have

$$u(x) = \sum_{k=1}^{N+L} \alpha_k h_k(x, x_k) u^* f(x_k) + \sum_{k=N+1}^{N+L} \beta_k p_k(x, x_k) u^* = \sum_{j=1}^{N+L} \gamma_k \pi_k(x, x_k), \qquad (13)$$

where $\alpha$, $\beta$ and $\gamma$ are unknown weighting coefficients, $h_k$ and $p_k$ are the indefinite functions. $hu^*f$, $pu^*$ and $\pi$ can be respectively interpreted as the radial basis functions for the DRM, boundary RBF methods such as the MFS and BKM, and domain-type RBF method [5]. Constructing RBF from numerical integral view of point is named as the general solution RBF (GS-RBF) [4]. How to choose function $h$ and $p$ will be left to discuss later.

RBF $\phi$ in Eq. (10) can be seen as $hu^*f$ in Eq. (13), which establishes a direct relationship between the RBF and Green integral. If $f(x)$ in Eq. (1) includes some point forcing terms, the RBF approximate representation can be stated as

$$u(x) = \sum_{j=1}^{N+L} \gamma_k \pi_k(x, x_k) + \sum_{l}^{q} Q_l u^*(x, x_l), \qquad (14)$$

where $q$ is the number of concentrated sources, $Q_l$ is the related magnitude.

### 2.2. QMC evaluation of particular solution

In contrast to the foregoing DRM strategy using the RBF approximation to Eq. (4), the Quasi-Monte Carlo scheme [3] directly evaluates domain integral of Eq. (11), which is in fact equivalent to the particular solution. The

QMC solution formula is given by

$$u_p(x) = \frac{V}{M}\sum_{k=1}^{M} u^*(x,x_k) f(x_k), \quad (15)$$

where $V$ is the volume of the domain, $M=N+L$ is the total number of nodes. The key in the QMC is proper node placement [2,3].

### 2.3. Equivalence of RBF and QMC

By comparing Eqs. (15), (13), and (10), it is easily established

$$\sum_{k=1}^{M} \alpha_k \psi_k(x,x_k) = \sum_{k=1}^{M} \alpha_k h_k u^*(x,x_k) f(x_k)$$
$$= \frac{V}{M}\sum_{k=1}^{M} u^*(x,x_k) f(x_k). \quad (16)$$

Eq. (16) shows that the RBF is closely related with the QMC inherently.

## 3. Further analysis and conjecture

The analysis of accuracy, convergence and stability of the RBF are very mathematically difficult [1]. On the other hand, the construction of efficient QMC algorithm is under intense study and still an open research topic. The underlying relation of both is mutually useful to get some insights into either method.

### 3.1. Constructing RBF and placing nodes

For inside-domain source points, the GS-RBF creates the RBFs

$$\varphi(x,x_k) = h_k(x,x_k) u^*(x,x_k) f(x_k). \quad (17)$$

Function $h$ assures that the RBF has enough differential continuity when $u^*$ is a singular fundamental solution. Singularity can also be eliminated if we place all response and source points differently or non-singular general solution is used [4]. In such cases, $h$ can be understood a RBF weight function. $h=r^{2m}$ is a convenient choose where $r$ is the Euclidean distance. One can find that the generalized thin plate splice (TPS) function $r^{2m}\ln(r)$ is a special case of the GS-RBF for the Laplace and biharmonic operators in 2D problems. $r^{2N+1}$ should be recommended for higher-dimensional problems.

Although Duchon [6] showed that the TPS with linear constraints is optimal for biharmonic operator interpolation, the linear constraints are not considered necessary if we modify TPS to $r^2\ln(r)+r^2+1$ according to the GS-RBF approach, where $\ln(r)$ general solution of biharmonic operator is not used due to its singularity. As of the boundary source points, we suggest

$$\varphi(x) = p_k(x,x_k) u^*(x,x_k) \quad (18)$$

as the RBF. $p$ functions like $h$ in Eq. (17).

What follows is an in-depth analysis of $h$ and $p$. The method of reduction of variance which includes importance sampling and stratified sampling are important to improve the accuracy of the MQC. Consider the integrand

$$\int w d\Omega = \int (w/g) g d\Omega = \int R g d\Omega, \quad (19)$$

where $gd\Omega$ denotes density which means the variable transformation method using the indefinite integral of $g$. The basic idea is to try to reduce function $R$ as close as possible to constant so that the accuracy of the QMC integral is considerably increased. This suggests that $h$ and $p$ should be positive definition with purpose averaging integrand function. Gaussian RBF may be interpreted as Gaussian probability distributions which are cases in some practical situations. Some sophisticated QMC techniques can be effectively used to construct efficient RBF within this framework. Also various singularity removal techniques in the BEM can be employed to construct RBF and QMC.

One immediate question is which of the QMC algorithm is a counterpart of the known multiquadratic (MQ) with spectral accuracy. For Green integral (11), we have

$$\int f(x) u^*(r,x) d\Omega = \int \sqrt{r^2+c^2} \frac{f(x) u^*(r,x)}{\sqrt{r^2+c^2}} d\Omega$$
$$= \int \sqrt{r^2+c^2} g(r,x,c) d\Omega, \quad (20)$$

Comparing to Eq. (19), the shape parameter $c$ of the MQ is here interpreted as giving the random points different weights corresponding to their position relative to other points in the QMC. Gamma distribution density may be the closest to the MQ. It is conceived that if the quasi-random nodes are weighted through varying parameter $c$ to keep the MQ constant as possible, the QMC algorithm can be the spectral convergence. However, this scheme may not be feasible due to the difficulty in indefinite integral of $g$. In contrast, the flexibility of the RBF lies in that it always can enforce arbitrary variable transformation in the reduction of variable method via collocation technique. This finding give some explanations to Kansa' variable shape parameter MQ [5].

The RBF may be applicable to evaluate any high-dimensional integral. Let integrand function in Eq. (19)

$$f(x) = w(x)/u^*(s_i, x) \qquad (21)$$

where $s_i$ is one specified node. Integrand (19) can be transformed to the evaluation of particular solution of Possion equation (1) with $f(x)$ defined in Eq. (21) as forcing term and zero Dirichlet boundary condition. Seeming inefficiency of this strategy due to the solution of a PDE can be eased if we use the spectral accuracy MQ as the RBF in the DRM. In addition, the multipole, multigrid and wavelet techniques can be used to reduce effort in the inversion of the RBF interpolation matrix to $O(M\log(M))$ [7].

Chen et al. [4] found that if we replace distance variable $r$ by $\sqrt{r^2+c^2}$ in the GS-RBF, we got the prewavelet RBF with the exponential accuracy. For example, the TPS is modified as $r^2 \ln\sqrt{r^2+c^2}$.

On the other hand, $hu^*f$ (or $pu^*$) from integral weight function viewpoint should be orthogonal with respect to node distribution just as in multidimensional orthogonal wavelet series. For simplicity, considering an ideal case where integrand function can be expressed only in terms of radial variable, we can reduce it to 1D problem

$$\int_a^b h(r)u^*(r)f(r)dr = \sum_{k=1}^{M} v_k h(r_k)u^*(r_k)f(r_k), \qquad (22)$$

where $v_k$ is weight coefficients. $hu^*$ is seen as the weight function in Gauss integration which can be chosen to remove integral singularities and make integral smooth. In the RBF case, the distribution of variable $r$ changes with response point. In this case we have abscissas at first and then we need to determine the weight function with a recurrence procedure to maximize the degree of accuracy with radical orthogonality. This is more or less similar to Gauss-Kronrod quadrature.

Fedoseyev et al [8] also found that compared with uniform nodes, the use of non-uniform ones are often much more efficient in the domain RBF solution of PDE. Furthermore, Fedoseyev et al. [9] addressed the severe accuracy drop of RBF PDE solutions at nodes neighboring boundary. Powell [7l] also found that the RBF geometry interpolation encounters evident accuracy loss around boundary. In fact, such edge effects are well known in polynomial interpolation and differential approximation, which are due to the geometry and system equation discontinuities at boundary. Remedy for geometry discontinuity is to place the nodes inclining boundary such as the use of zeros spacing of the Lobatto, Chebyshev or legendre polynomials. This methodology also works for the system discontinuity problem. Therefore, the domain RBF scheme for PDE should use non-uniform nodes. A quantitative measure for node distributions is

$$\sigma(x_i) = \frac{\Omega}{M} \sum_{k=1}^{M} \|x_i - x_k\|. \qquad (23)$$

Geometry discontinuity can be understood that $\sigma$ value of distinct node is different, especially between boundary and central nodes. Fedoseyev et al. [9] seem to solve the system discontinuity problem very well by collocating the governing equation at boundary points.

Similar to the RBF, it is also well known that the node placement has crucial effect on efficiency and accuracy of the QMC. The RBF and QMC can be mutually compared to place nodes optimally, especially for complex high-dimension domain. The edge discontinuous effects strongly suggest that the density of node distribution in the QMC should be inclined to boundary rather than evenly spaced. Also the more nodes should be placed in regions of highly gradient integrand function.

**3.2. Conjecture on the RBF error bounds**

The RBF error bounds available now do not include dimension effect. Kansa and Hon [5] observed that the RBF is dimension independent numerically. They stressed that the convergence order of the RBF increases as the spatial dimension increases, which is featured by error bound $O(M^{(d+1)})$. Here $d$ means the spatial dimension. However, they did not give any evidence supporting this error estimate. Our numerical experiments did not justify this error formula.

It is very interesting to note that there exist the same error behaviors between some RBFs and QMC or MC algorithms. For example, error bounds for the linear RBF and the classical Monte Carlo method are the same $O(M^{-1/2})$ which $M$ is the number of nodes, while error bounds for the TPS RBF and QMC in 2D problems are the same $O(M^{-1}(\log M))$. It is well known that the QMC has error bound

$$err = O\left(M^{-1}(\log M)^{d-1}\right) \qquad (24)$$

for $d$-dimension problems. The following error bound conjecture of the RBF is intuitively proposed.

**Conjecture 1**: The RBF error bounds can be characterized by $O(M^{-\eta}(\log M)^{d-1})$, where $d$ is dimension, $\eta$ is 1/2 for linear RBF, 1 for the TPS in the 2D problems, and $n$ for the MQ.

The dimension affect on the RBF accuracy is here featured by $O(\log M)^{d-1}$ rather than $O(M^{(d+1)})$ in [5]. The notorious dimension curse in the other numerical techniques can be characterized by

$$err = O\left(M^{-\kappa/d}\right), \qquad (25)$$

where $\kappa$ denotes the accuracy order of algorithm. By comparing formulas (24) and (25), one can easily conclude that the QMC method has very visible advantages in accuracy for higher dimensional problems. For example, if the accuracy of some FEM or FDM scheme is three-order. Then, the scheme has the same efficiency as the QMC for 3D problems, while the method will be much less efficient for higher dimension problems. In addition, the time-space RBF [4]

$$r_j = \sqrt{(x-x_j)^2 + c(t-t_j)^2}, \qquad (26)$$

can be applied to the wave propagation problem, where $c$ is the wave velocity. The GS-RBF using time-dependent fundamental or general solution also eliminates time dependence of the convection-diffusion problems. In this way, the dimensional accuracy advantage of the RBF scheme is extended to time dimension.

**4. Kernel function and RBF**

Chen et al. [4] pointed out that the RBF has very close relationship with the kernel function of integral equation based on an integral analysis of the RBF solution of PDE. The RBF is now also widely employed in network, data processing, and inverse problems. These problems can be expressed by various integral equations. For instance,

$$Rh = \int_D R(x,y)h(y)dy = f(x), x \in \overline{D} \subset IR^n \qquad (27)$$

is the basic equation of stochastic optimization theory, where the $R(x,y)$ is the kernel of a positive rational function of an elliptic self-adjoint operator $L$ on $L^2(IR^n)$ [10]. The creating RBF for this problem should involve the kernel function $R$ as done in Eq. (27).

Another instance is geometry generation which can be understood an elliptic boundary value problems relating to the Laplace and biharmonic operators [11]. The fact that the $\ln(r)$ in TPS is an essential ingredient of the fundamental solutions of 2D Laplace and biharmonic operator explains why the TPS is widely used in 2D geometry generation problems.

According to the above discussions, the brutal use of the TPS and MQ should be discouraged. The kernel function is found important to construct efficient operator-dependent RBFs.

## 5. Concluding remarks

It is noted that $r^{2N}$ are usually not efficient compared with $r^{2N+1}$ as the RBF. The formal is a polynomial interpolation, while latter as well as the TPS, MQ, Gauss, and GS-RBF are all irrational interpolation technique. Therefore, the interpolation analysis of RBF is not an easy task. Instead, this note shows that the QMC may provide an alternative promising research approach for the RBF. Some perplexing issues pertaining to RBF approach such as the accuracy, stability, knot generation, and convergence may be handled by means of the existing results relating to the QMC. On the other hand, some RBF-related methods such as the multipole, quasi-interpolation and multilevel technique with compactly-supported RBF may be extended to the QMC. The research of the RBF and QMC can converge in some ways in the near future.